\documentclass{amsart}
\usepackage[final]
{showkeys}
\usepackage{graphicx}
\usepackage{enumerate}
\usepackage{url}
\usepackage{hyperref}

\theoremstyle{theorem}
\newtheorem{theorem}{Theorem}

\newtheorem{corollary}[theorem]{Corollary}

\theoremstyle{definition}

\theoremstyle{definition} 
\newtheorem*{remark*}{Remark}

\newcommand{\E}{\operatorname{\mathsf{E}}}
\renewcommand{\P}{\operatorname{\mathsf{P}}}

\newcommand{\A}{\mathcal{A}}
\newcommand{\B}{\mathcal{B}}

\newcommand{\Om}{\Omega}

\begin{document}


\title{Exact lower bound on an ``exactly one'' probability}
\author{Iosif Pinelis}
\address{Department of Mathematical Sciences, Michigan Technological University, Houghton, Michigan 49931}
\email{ipinelis@mtu.edu}






\begin{abstract}
The exact lower bound on the probability of the occurrence of exactly one of $n$ random events each of probability $p$ is obtained.  
\end{abstract}


\maketitle

%




\section{Introduction, summary, and discussion}\label{intro}
Suppose $A_1,\dots,A_n$ are random events each of probability $p$. Let $E$ denote the event that exactly one of the events $A_1,\dots,A_n$ occurs. 

If the $A_i$'s are independent then, by the binomial probability mass function formula (see e.g.\ \cite[Section~1.3]{borovkov-prob}), $\P(E)=npq^{n-1}$, where $q:=1-p$. So, in the ``independent'' case, $\P(E)$ attains its maximum, equal $(1-1/n)^{n-1}\big[\underset{n\to\infty}\longrightarrow1/e\big]$, at $p=1/n$.  

What will happen with $\P(E)$ when the $A_i$'s are only assumed to be pairwise independent? Of course, already for $n=3$, the pairwise independence of the $A_i$'s does not imply their ``complete'' independence. Feller \cite[page~126]{feller_vol1} wrote: ``Actually such occurrences [of pairwise independence but not ``complete'' independence] are so rare that their possibility passed unnoticed until S.\ Bernstein constructed an artificial example. It still takes some search to find a plausible natural example.'' 
This is followed (\cite[page~127]{feller_vol1}) by an example 
of three pairwise independent events that are not ``completely'' independent. 
Another such example, \cite[Example~2.3.3]{borovkov-prob} -- ascribed in \cite{borovkov-prob} to Bernstein, actually appears more common and natural than the mentioned example on page~127 in \cite{feller_vol1}. 

One may want to dispute the assertion that occurrences of pairwise independence without ``complete'' independence  are rare. Indeed, the definition of the independence of three events $A,B,C$ consists of the following four equations: $\P(A\cap B)=\P(A)\P(B)$, $\P(B\cap C)=\P(B)\P(C)$, $\P(A\cap C)=\P(A)\P(C)$, and 
$\P(A\cap B\cap C)=\P(A)\P(B)\P(C)$.  
The first three of these four equations define the pairwise independence. The probabilities of the events $A,B,C$ and of their pairwise and triple intersections can all be expressed as the sums of the probabilities of certain pieces of the partition of the sample space (say $\Om$) generated by the events $A,B,C$. There are $2^3=8$ pieces of this partition, with $8$ corresponding probabilities, which may be considered as nonnegative real variables tied just by one equation -- stating that the sum of these $8$ probabilities is $1$. Thus, we have $4+1=5$ equations with $8$ unknowns, which leaves us $8-5=3$ degrees of freedom, which one can easily use to show that none of the four equations defining the independence of the events $A,B,C$ may be dropped without altering the notion of independence. In particular, this way it is easy to see that the pairwise independence does not imply the ``complete'' independence. Moreover, now it seems plausible that -- in the case of three events $A,B,C$ -- the dimension of the semi-algebraic set \cite{bochnak-etal} (in the mentioned $8$ variables) corresponding to the ``complete'' independence is less by $1$ than the dimension of the semi-algebraic set corresponding to the pairwise independence. More generally, for any natural number $n$ of events $A_1,\dots,A_n$, the difference between the corresponding dimensions appears to be $
2^n-1-n-n(n-1)/2\sim2^n$ (as $n\to\infty$). Here, $2^n$ appears as the number of all equations of the form $\P\big(\bigcap_{j\in J}A_j\big)=\prod_{j\in J}\P(A_j)$
%
for $J\subseteq[n]:=\{1,\dots,n\}$. 
Of these $2^n$ equations, $1+n$ equations -- for the sets $J\subseteq[n]$ of cardinalities $0$ and $1$ -- are trivial; and $n(n-1)/2$ of the $2^n$ equations define the pairwise independence of the $n$ events. Thus, there are $
2^n-1-n-n(n-1)/2$ ``nontrivial'' equations defining the independence of the $n$ events, in addition to the $n(n-1)/2$ ``nontrivial'' equations defining the pairwise independence. 

From the just described viewpoint, the occurrences of ``complete'' independence constitute an infinitesimally thin slice among the occurrences of pairwise independence. 
Therefore, it may seem very surprising that the strong law of large numbers (SLLN) for identically distributed random variables with a finite mean turns out to hold assuming only pairwise independence, as was 
demonstrated comparatively very recently by Tao \cite[Remark 2]{tao-SLLN}.

In this note it will be shown that, in contrast with the latter SLLN result, the ``exactly one'' probability $\P(E)$ may be quite sensitive to the distinction between the pairwise independence and the ``complete'' independence: 

\begin{theorem}\label{th:} 
For each natural $n$ and each $p\in[0,1]$,
\begin{equation}\label{eq:}
	\min\P(E)=P_{n,p}:=np\big(1-(n-1)p\big)_+,
\end{equation}
where the minimum is taken over all pairwise independent events $A_1,\dots,A_n$ each of probability $p$, and $x_+:=\max(0,x)$ for real $x$. 
\end{theorem} 

We see that, in contrast with the ``completely independent case'', for just pairwise independent events $A_1,\dots,A_n$ the probability $\P(E)$ can be $0$ for any $n\ge2$ and any $p\ge1/(n-1)$. 
If we consider the special value $1/n$ of $p$ -- at which, as was noted, $\P(E)$ attains its \emph{maximum} value $(1-1/n)^{n-1}\approx1/e$ in the ``completely independent'' case --   
then for just pairwise independent events $A_1,\dots,A_n$ we have $\min \P(E)=1/n\to0$. 
However, if e.g.\ $p=c/n$ with a fixed $c\in(0,1)$, then in both cases -- of the ``complete'' independence and of the pairwise independence -- the probability $\P(E)$ stays away from $0$. 
So, $\P(E)$ will necessarily be of the same order of magnitude (for large $n$) in both cases only if $p$ small -- more specifically, if $p$ stays below $c/n$ for some fixed $c\in(0,1)$. 

\begin{figure}[htbp]
	\centering
		\includegraphics[width=1.00\textwidth]{
		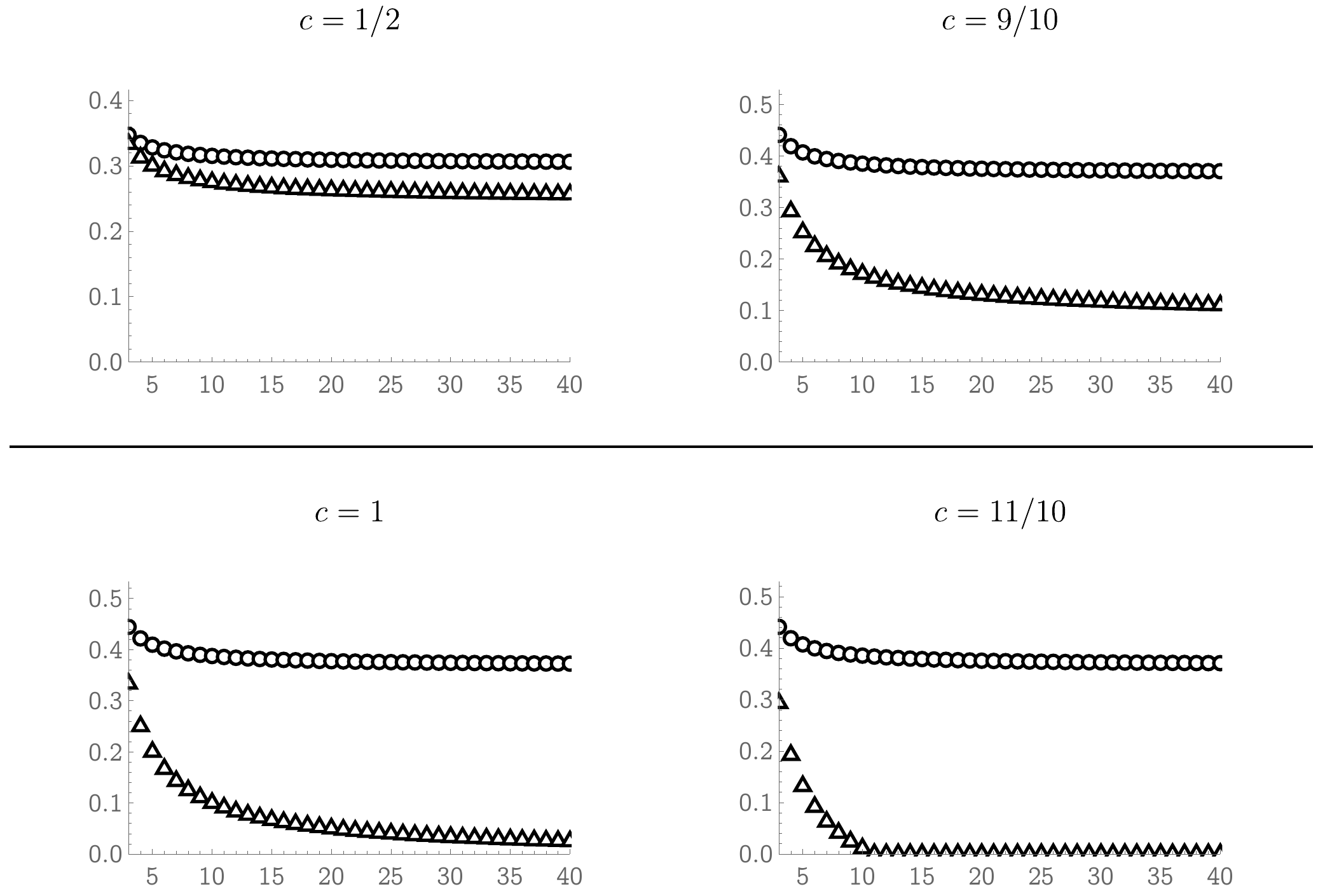}
	\caption{Graphs of the values of $\P(E)$ for $p=c/n$}
	\label{fig:pic}
\end{figure}

This is illustrated in Fig.\ \ref{fig:pic}, which shows the graphs of the values of $\P(E)$ (the vertical axis) in the ``completely independent case'' (circles) and in the ``pairwise independent case'' (triangles) for $n\in\{3,\dots,40\}$ (the horizontal axis), $p=c/n$, and $c\in\{1/2,\,9/10,\,1,\,11/10\}$. 

\section{Proof of Theorem~1
}\label{proof}
For $n=1$, Theorem~\ref{th:} is trivial. So, in what follows assume $n\ge2$. 

For each $j\in[n]$, let 
\begin{equation*}
X_j:=1_{A_j},	
\end{equation*}
the indicator of the event $A_j$. Let 
\begin{equation*}
N:=X_1+\dots+X_n, 	
\end{equation*}
the number of the events $A_1,\dots,A_n$ that occurred. Then 
\begin{equation}\label{eq:E=}
	E=\{N=1\}. 
\end{equation}
Note that $\E X_j=p$ and (by the pairwise independence) $\E X_jX_k=p^2+pq\,1(j=k)$ for all $j$ and $k$ in $[n]$. 
Now we have a perhaps unexpected use of 
%
the Chebyshev--Markov inequality (see e.g.\ \cite[Theorem~4.7.4]{borovkov-prob}): 
\begin{align*}
	\P(N\ne1)&=\P((N-1)^2\ge1) \\ 
&\le\E(N-1)^2 \\ 
&=\E N^2-2\E N+1\\ 
&=\sum_{j,k\in[n]}\E X_jX_k-2\sum_{j\in[n]}\E X_j+1 \\ 
&=n^2p^2+npq-2np+1 \\ 
&=1-np\big(1-(n-1)p\big). 
\end{align*}
Therefore and because $\P(N=1)\ge0$, we see that  
\begin{equation*}
	\P(N=1)\ge np\big(1-(n-1)p\big)_+=P_{n,p}; 
\end{equation*}
So, in view of \eqref{eq:E=}, $P_{n,p}$ is a lower bound on $\P(E)$; 
cf.\ \eqref{eq:}. 

It remains to show that this 
lower bound 
is attained, for each natural $n\ge2$ and each $p\in[0,1]$. 
To do this, introduce the events 
\begin{equation*}
	C_J:=\Big(\bigcap_{j\in J}A_j\Big)\cap\Big(\bigcap_{j\in[n]\setminus J}(\Om\setminus A_j)\Big)
\end{equation*}
for $J\subseteq[n]$. These events constitute a partition of the sample space $\Om$.  
Moreover, for each $m\in\{0\}\cup[n]$, 
\begin{equation}\label{eq:N=m}
	\{N=m\}=\bigcup_{\substack{J\subseteq[n],\\ |J|=m}}C_J,
\end{equation}
where $|J|$ denotes the cardinality of the set $J$. Also, 
\begin{equation}\label{eq:A1,A2}
	A_1=\bigcup_{\substack{J\subseteq[n],\\ J\supseteq\{1\}}}C_J \quad\text{and}\quad  
	A_1\cap A_2=\bigcup_{\substack{J\subseteq[n],\\ J\supseteq\{1,2\}}}C_J. 
\end{equation}

For each $m\in\{0\}\cup[n]$, let us assign the same probability, say $x_m$,
to each event $C_J$ with $J\subseteq[n]$ such that $|J|=m$; then, by \eqref{eq:N=m}, 
\begin{equation}\label{eq:P(N=m)}
	\P(N=m)=\binom nm x_m.
\end{equation}
So, there will exist a probability space supporting such an assignment of probabilities to the $C_J$'s 
%
if and only if $x_m\ge0$ for all $m\in\{0\}\cup[n]$ and 
\begin{equation}\label{eq:sum=1}
	\sum_{m=0}^n\binom nm x_m=1;
\end{equation}
this follows because the set of values of the random variable $N$ is the set $\{0\}\cup[n]$.  


Then, in view of \eqref{eq:A1,A2}, we also have 
\begin{equation*}
	\P(A_1)=\sum_{m=1}^n \sum_{\substack{J\subseteq[n],\\ J\supseteq\{1\},\\ |J|=m}}\P(C_J)
	=\sum_{m=1}^n \binom{n-1}{m-1} x_m 
\end{equation*}
(which is actually the value of $\P(A_j)$ for all $j\in[n]$) and 
\begin{equation*}
	\P(A_1\cap A_2)=\sum_{m=1}^n \sum_{\substack{J\subseteq[n],\\ J\supseteq\{1,2\},\\ |J|=m}}\P(C_J)
	=\sum_{m=1}^n \binom{n-2}{m-2} x_m 
\end{equation*}
(which is actually the value of $\P(A_i\cap A_j)$ for all distinct $i$ and $j$ in the set $[n]$). 
Now the conditions that $\P(A_j)=p$ for all $j\in[n]$ and the $A_j$'s are pairwise independent can be rewritten as 
\begin{equation}\label{eq:=p,p^2}
	\sum_{m=1}^n \binom{n-1}{m-1} x_m=p\quad\text{and}\quad
	\sum_{m=1}^n \binom{n-2}{m-2} x_m=p^2. 
\end{equation}

Now take any $p\in[0,1]$. Then there is some $k\in[n-1]$ such that 
\begin{equation}\label{eq:p case}
	\frac{k-1}{n-1}\le p\le\frac{k}{n-1}. 
\end{equation}
For such a number $k\in[n-1]$, 
let 
\begin{equation}\label{eq:x_m:=}
	x_m:=\left\{
	\begin{alignedat}{2}
	&\frac{np}k\,\big(k-(n-1)p\big)\Big/\binom nk &&\text{\quad if }m=k, \\ 
	&\frac{np}{k+1}\,\big((n-1)p-(k-1)\big)\Big/\binom n{k+1} &&\text{\quad if }m=k+1, \\ 
	&0&&\text{\quad if }m\in[n]\setminus\{k,k+1\}. 
	\end{alignedat}
	\right.
\end{equation}
Then, in view of condition \eqref{eq:p case}, $x_m\ge0$ for all $m\in
[n]$. Also, then straightforward calculations show that conditions \eqref{eq:=p,p^2} hold and 
\begin{equation}\label{eq:s:=}
	s:=\sum_{m=1}^n\binom nm x_m=\frac{n p\big(2 k-(n-1)p\big)}{k (k+1)}\le1. 
\end{equation}
(The latter inequality 
is elementary. To prove it, one may first note that the maximum in $p$ of the ratio in \eqref{eq:s:=} is $\frac{k n}{(k+1) (n-1)}$, which increases in $k\in[n-1]$ to $1$.) 
Therefore, one can satisfy condition \eqref{eq:sum=1} by letting $x_0:=1-s\ge0$, so that the condition $x_m\ge0$ for all $m\in\{0\}\cup[n]$ holds as well.

Furthermore, it follows from \eqref{eq:E=}, \eqref{eq:P(N=m)}, \eqref{eq:p case}, \eqref{eq:x_m:=}, and the definition of $P_{n,p}$ in \eqref{eq:} that 
\begin{align*}
	\P(E)=\P(N=1)=n x_1
&	=\left\{
	\begin{alignedat}{2}
	&np\big(1-(n-1)p\big)&&\text{\quad if }0\le p\le\frac1{n-1}, \\ 
	&0&&\text{\quad otherwise}. 
	\end{alignedat}
	\right. \\
&=np\big(1-(n-1)p\big)_+=P_{n,p}. 	
\end{align*} 

This shows that the lower bound $P_{n,p}$ on $\P(E)$ is indeed attained, which completes the proof of Theorem~\ref{th:}. \qed

We have the following easy corollary of Theorem~\ref{th:}: 

\begin{corollary}\label{cor:}
In the conditions of Theorem~\ref{th:}, the best lower bound on $\P(N=n-1)$ is $P_{n,q}$ (cf.\ \eqref{eq:E=}). 
\end{corollary} 

To see why this corollary holds, switch from the ``successes'' $A_j$ to the ``failures'' $\Om\setminus A_j$, and also interchange the roles of $p$ and $q=1-p$. 

\bigskip

There are a number of further questions that one may ask concerning Theorem~\ref{th:}, including the following: 

\begin{enumerate}[1.]
	\item Assuming still that $A_1,\dots,A_n$ are pairwise independent events each of probability $p$, what is the best \emph{upper} bound on $\P(E)=\P(N=1)$? More generally, for each $m\in\{0\}\cup[n]$, under the same conditions on the $A_j$'s, what are the best lower and upper bounds on $\P(N=m)$?
	\item The same questions as above, but assuming, more generally, that the $A_j$'s are $r$-independent for some $r\in\{2,\dots,n-1\}$, i.e., assuming that for any $J\subseteq[n]$ with $|J|=r$ the family $(A_j)_{j\in J}$ is independent.  
	\item The same questions as above, but assuming, more generally, that the probabilities $\P(A_j)$ have possibly different prescribed values $p_j$, for $j\in[n]$. 
	\item Yet more generally, let $\B$ be any subset of the algebra (say $\A$) generated by events $A_1,\dots,A_n$. Suppose that the probabilities $\P(B)$ have prescribed values, say $p_B$, for all $B\in\B$. Take any $A\in\A$. What are the best lower and upper bounds on $\P(A)$ in terms of the $p_B$'s?
\end{enumerate}

Looking back at the proof of Theorem~\ref{th:} and recalling the discussion in Section~\ref{intro}, one can see that 
%
all the further problems listed above are ones of linear programming in a space of dimension exponentially growing with $n$, with the values of the $\P(C_J)$'s for $J\subseteq[n]$ as the variables. Therefore and because the above proof of Theorem~\ref{th:}, with all its parts fitting together quite tightly, already 
was
not easy to 
devise, all these problems seem hard to tackle theoretically or even computationally. 

%


\section{Conclusion}
As we saw in Section~\ref{intro}, the condition of the ``complete'' independence of $n$ events, oftentimes assumed quite casually, actually involves $\sim2^n$ equations, which are practically impossible to test well even for rather moderate values of $n$, such as $n=40$. In contrast, the pairwise independence of $n$ events involves only $n(n-1)/2$ conditions. It may therefore be of value and interest to know how much the consequences of these two kinds of independence may differ from each other in various settings. It was noted in Section~\ref{intro} that, at least as far as the most common version of the strong law of large numbers (for identically distributed random variables with a finite mean) is concerned, the pairwise independence is just as good as the ``complete'' independence of the random variables. 
In stark contrast with that, the ``exactly one'' probability may be quite sensitive to the distinction between the pairwise independence and the ``complete'' independence, as shown in this note. 

It is hoped that this small study may stimulate further research into the other aspects of the difference between the ``complete'' independence and, on the other hand, the pairwise independence (or, more generally, the $r$-independence for some $r\in\{2,\dots,n-1\}$). Also, perhaps some of the further questions enumerated at the end of Section~\ref{proof} will attract attention of other researchers. Finally, the methods presented in this note might turn out to be of use in other optimization problems in probability, statistics, and perhaps elsewhere, especially where the ``complete'' independence is in doubt.

\bibliographystyle{amsplain}

\bibliography{P:/pCloudSync/mtu_pCloud_02-02-17/bib_files/citations01-09-20}

\end{document}